\documentclass[12pt, reqno]{amsart}

\usepackage{amssymb}  
\usepackage{array}
\usepackage{enumitem}
\usepackage{url}
\usepackage{verbatim}

\theoremstyle{plain}
\newtheorem{theorem}{Theorem}[section]
\newtheorem{lemma}[theorem]{Lemma}

\newtheorem{proposition}[theorem]{Proposition}
\newtheorem{corollary}[theorem]{Corollary}

\newtheorem{problem}[theorem]{Problem}

\theoremstyle{remark}

\numberwithin{equation}{section}

\newcommand{\seclabel}[1]{\label{sec:#1}}   
\newcommand{\thmlabel}[1]{\label{thm:#1}}   
\newcommand{\lemlabel}[1]{\label{lem:#1}}   
\newcommand{\corlabel}[1]{\label{cor:#1}}   
\newcommand{\prplabel}[1]{\label{prp:#1}}   
\newcommand{\prblabel}[1]{\label{prb:#1}}   
\newcommand{\eqnlabel}[1]{\label{eqn:#1}}   

\newcommand{\secref}[1]{\ref{sec:#1}}   
\newcommand{\thmref}[1]{\ref{thm:#1}}   
\newcommand{\lemref}[1]{\ref{lem:#1}}   
\newcommand{\corref}[1]{\ref{cor:#1}}   
\newcommand{\prpref}[1]{\ref{prp:#1}}   
\newcommand{\prbref}[1]{\ref{prb:#1}}   
\newcommand{\eqnref}[1]{\eqref{eqn:#1}} 

\newcommand{\by}[1]{\overset{\eqnref{#1}}=}  

\newcommand{\Aut}{\mathrm{Aut}}			
\newcommand{\Mlt}{\mathrm{Mlt}}         
\newcommand{\Inn}{\mathrm{Inn}}			
\newcommand{\Fix}{\mathrm{Fix}}         

\newcommand{\setof}[2]{\{#1\,|\,#2\}}   

\newcommand{\ld}{\backslash}					
\newcommand{\inv}{^{-1}}						
\newcommand{\sbl}[1]{\langle#1\rangle}			
\newcommand{\id}[1]{\mathrm{id}_{#1}}           
\newcommand{\Al}{\textsc{A}}                    
\newcommand{\AIP}{\textsc{AIP}}                 
\newcommand{\nlp}{\diamond}                     

\title{The structure of commutative automorphic loops}

\author{P\v{r}emysl Jedli\v{c}ka}
\author{Michael Kinyon}
\author{Petr Vojt\v{e}chovsk\'y}

\address[Jedli\v{c}ka]{Department of Mathematics \\
Faculty of Engineering \\ Czech University of Life Sciences \\
Kam\'yck\'a 129 \\ 165 21 Prague 6–-Suchdol \\ Czech Republic}
\email[Jedli\v{c}ka]{\url{jedlickap@tf.czu.cz}}

\address[Kinyon and Vojt\v{e}chovsk\'y]{Department of Mathematics \\
University of Denver \\ 2360 S Gaylord St \\ Denver, Colorado 80208 USA}
\email[Kinyon]{\url{mkinyon@math.du.edu}}
\email[Vojt\v{e}chovsk\'y]{\url{petr@math.du.edu}}

\begin{document}

\begin{abstract}
An \emph{automorphic loop} (or \emph{A-loop}) is a loop whose inner mappings
are automorphisms. Every element of a commutative A-loop generates a group, and
$(xy)^{-1} = x^{-1}y^{-1}$ holds. Let $Q$ be a finite commutative A-loop and
$p$ a prime. The loop $Q$ has order a power of $p$ if and only if every element
of $Q$ has order a power of $p$. The loop $Q$ decomposes as a direct product of
a loop of odd order and a loop of order a power of $2$. If $Q$ is of odd order,
it is solvable. If $A$ is a subloop of $Q$ then $|A|$ divides $|Q|$. If $p$
divides $|Q|$ then $Q$ contains an element of order $p$. If there is a finite simple
nonassociative commutative A-loop, it is of exponent $2$.
\end{abstract}

\maketitle

\section{Introduction}
\seclabel{intro}

A \emph{loop} $(Q,\cdot)$ is a set $Q$ with a binary operation
$\cdot$ such that (i) for each $x\in Q$, the
\emph{left translation}
$L_x : Q\to Q ; y\mapsto yL_x = xy$
and the \emph{right translation}
$R_x : Q\to Q ; y\mapsto yR_x = yx$
are bijections, and (ii) there exists $1\in Q$ satisfying
$1\cdot x = x\cdot 1 = x$ for all $x\in Q$.
The left and right translations generate the \emph{multiplication group}
$\Mlt(Q) = \sbl{L_x, R_x \mid x\in Q}$. The \emph{inner mapping group}
$\Inn(Q) = \Mlt(Q)_1$ is the stabilizer of $1\in Q$.
Standard references for the theory of loops are
\cite{Belousov, Bruck, Pflugfelder}.

A loop $Q$ is an \emph{automorphic loop} (or \emph{A-loop}) if every inner
mapping of $Q$ is an automorphism of $Q$, that is, $\Inn(Q) \leq \Aut(Q)$. Thus
the class of A-loops, which is certainly not the class of all loops, includes,
for instance, groups and commutative Moufang loops \cite{Bruck}.

The study of A-loops was initiated by Bruck and Paige \cite{BP}.
They obtained many basic structural results for A-loops and
also described some constructions.
The bulk of \cite{BP} was devoted to the (implicitly stated)
problem of whether every
\emph{diassociative} A-loop, that is, an A-loop in which every
$2$-generated subloop is a group, is a Moufang loop.
Affirmative answers were given by Osborn \cite{Osborn} in
the commutative case, and Kinyon, Kunen and Phillips \cite{KKP}
in the general case. Moufang A-loops have been used to
characterize a certain class of quasigroups \cite{Kepka},
and have been shown to have an affirmative answer for the
restricted Burnside problem \cite{PS}.

By contrast, the study of other classes of A-loops has lain quite fallow. In
this paper, we give a detailed structure theory for \emph{commutative} A-loops.
Here is a summary of our main results:

In \S\secref{prelims}, we present preliminary results which will be used
throughout the rest of the paper. Some of these results, such as the
power-associativity of commutative A-loops (Lemma \lemref{pa}) are already
known for arbitrary A-loops \cite{BP}, but we give short proofs to make the
present paper self-contained. Other results, such as the automorphic inverse
property (Lemma \lemref{aaip}) are new.

In \S\secref{odd-order}, we study commutative A-loops of odd order, i.e. finite
A-loops in which every element has odd order (Lemma \lemref{odd-2div}). The
multiplication group of a commutative A-loop contains a natural (but not at all
obvious) twisted subgroup (Lemma \lemref{twisted}). In the odd order case, this
enables us to construct a new loop operation on a commutative A-loop with the
property that powers in the new loop coincide with powers in the original loop
(Lemma \lemref{bruck}). The new loop is in fact a \emph{Bruck loop}, and we
exploit this fact to establish \emph{Lagrange} and
\emph{Cauchy} theorems for commutative A-loops of odd order (Propositions
\prpref{lagrange} and \prpref{HSC}). Our main result in \S\secref{odd-order} is
the \emph{Odd Order Theorem}: every commutative A-loop of odd order is solvable
(Theorem \thmref{ft}).

In \S\secref{squares}, we turn to a property trivially satisfied in abelian
groups and valid in commutative Moufang loops thanks to dissociativity: the
product of squares is a square. This turns out be true in commutative A-loops
as well (Theorem \thmref{SquareRoot}), despite the fact that the naive formula
$x^2 y^2 = (xy)^2$ does not hold in general. Instead, $x^2 y^2 = (x\nlp y)^2$
for a rather complicated binary operation $\nlp$; in the Moufang case, $\nlp$
coincides with the original operation. Following the same philosophy as in the
odd order case, we study the new operation $\nlp$ and note that it defines a
commutative, power-associative loop on the same underlying set as the original
commutative A-loop. In the odd order case, $\nlp$ yields an isomorphic copy of
the original loop (Lemma \lemref{2div-iso}), but at the other extreme where
every element has order a power of $2$, the new loop operation $\nlp$ turns out
to have strong structural properties, as we will show in later sections.

In \S\secref{decomposition}, we prove a \emph{Decomposition Theorem}: every
finite commutative A-loop is a direct product of a subloop of odd order and a
subloop in which every element has order a power of $2$ (Theorem
\thmref{Decomposition}). This is a generalization of the familiar decomposition
theorems in abelian groups and commutative Moufang loops. Unlike in those
cases, however, no further decomposition is possible: commutative A-loops of
odd order are not necessarily direct products of $p$-loops for odd $p$.

In \S\secref{Exp2}, we examine commutative A-loops of exponent $2$. This
special case is of particular importance because of a straightforward
consequence of the Decomposition Theorem and the Odd Order Theorem, namely that
a finite, simple, commutative A-loop is either a cyclic group of odd prime
order or it has exponent $2$ (Proposition \prpref{simple}). To study the
exponent $2$ case, we return to the new loop operation $\nlp$ introduced in
\S\secref{squares}, and prove the main result of \S\secref{Exp2}: if $Q$ is a
finite, commutative A-loop of exponent $2$, then $(Q,\nlp)$ is an elementary
abelian $2$-group (Theorem \thmref{abelian-group}). An immediate corollary of
this is that a commutative A-loop of exponent $2$ has order a power of $2$
(Corollary \corref{exp2-order}).

In \S\secref{p-loops}, we briefly examine $p$-loops. The main result is that
the two reasonable definitions of this notion coincide for commutative A-loops,
that is, a finite commutative A-loop has order a power of $p$ if and only if
every element has order a power of $p$ (Theorem \thmref{p-loop}). For $p$ odd,
this is a consequence of the Lagrange and Cauchy theorems. For $p = 2$, it
follows from the Decomposition Theorem and the fact that it has already been
observed in the exponent $2$ case. We now easily derive the \emph{Lagrange}
and \emph{Cauchy Theorems} for all finite commutative
A-loops (Theorem \thmref{LHSC}).

Finally, in \S\secref{open} we state three open problems. The first, which we
expect to generate a great deal of interest in loop theory, is whether there
exists a nonassociative, finite simple commutative A-loop (Problem
\prbref{simple}). The results in this paper already tell us a great deal about
the structure such a loop must have. The second problem (Problem
\prbref{center}) is whether every commutative A-loop of odd prime power order
has a nontrivial center, that is, whether the loop is centrally nilpotent.
The third (Problem \prbref{hall}) addresses the existence of Hall and
Sylow subloops.

We should note that the variety of commutative A-loops is vast compared to the
variety of abelian groups. There exist many nonassociative examples even under
very restrictive conditions, such as in the case of commutative A-loops of
exponent two. While every A-loop of prime order $p$ is isomorphic to the cyclic
group of order $p$, a class of nonassociative commutative A-loops of order $pq$
($2<p<q$ primes) was found by Dr\'apal \cite{Drapal}. A survey of known
constructions and the classification of commutative A-loops of small orders
will appear in the planned sequel \cite{JKV2} to this paper. In \cite{JKV2}, we
also give an example of a commutative A-loop of order $16$ that is not
centrally nilpotent.

The main idea of this paper is to associate a new loop operation with the
original loop. In the odd order case, where the original loop is uniquely
$2$-divisible, this is a familiar approach \cite{Glauberman2}, \cite{FKP}.
However, in all earlier instances it was somewhat transparent what the
associated loop operation should be, unlike here. A common feature is to take
advantage of the unique square roots. We do not have access to square roots in
$2$-loops, but if for every $x$, $y$ there is $z$ such that $x^2y^2=z^2$
(Theorem \thmref{SquareRoot}), our novel idea is to declare $z$ to be a new
product of $x$ and $y$. As demonstrated in this paper, this approach is most
fruitful in case of commutative A-loops. Moreover, we now have some anecdotal
evidence that the connection is more profound, and that binary operations
associated in this or similar manner are deserving of a systematic
investigation in other varieties of loops.

The well-behaved structure theory of commutative A-loops belies the rather
technical lemmas on which it is based. Most of these lemmas involve detailed
equational reasoning, often obtained with the assistance of the automated
theorem prover Prover9 \cite{McCune}.

Finally, we should mention that many of our structural results for commutative A-loops of
odd order can be generalized to the noncommutative case. These generalizations
will appear elsewhere \cite{KKP2}.

\subsection{Notation}

Throughout the paper, let $Q$ denote a commutative loop with multiplication denoted by
juxtaposition and with neutral element $1$. Since all left translations are
bijections of $Q$, it is convenient to define the associated left division operation
by
\[
x\ld y = yL_x\inv
\]
for all $x,y\in Q$. It will also be useful to introduce the \emph{division permutations}
$D_x : Q\to Q$, $x\in Q$, defined by
\[
y D_x = y \ld x = x L_y\inv
\]
for all $x,y\in Q$. Note that $D_x^2 = \id{Q}$ for all $x\in Q$.
We will use the usual notation $x\inv = x \ld 1$ for the inverse of $x$, and we will
also use the \emph{inversion permutation} $J : Q\to Q$ defined by
\[
xJ = xD_1 = x\inv
\]
for all $x\in Q$.

To avoid excessive parenthesization, we will use the following convention.
The multiplication operation $\cdot$ will be less binding than left division,
which is, in turn, less binding than juxtaposition. For example, with this
convention, $ab\ld cd \cdot g\ld ef$ is unambiguously read as
$((ab)\ld (cd))(g\ld (ef))$. On the other hand, we shall certainly use
parentheses, brackets, \emph{etc.}, whenever they help to clarify an expression.

It is well known \cite{Bruck} that for commutative loops, the inner mapping
group $\Inn(Q)$ has a distinguished set of generators
\[
L_{x,y} = L_x L_y L_{yx}\inv
\]
for $x,y\in Q$. Using these generators, the A-loop condition can be expressed
as follows:
\[
(uv) L_{x,y} = uL_{x,y} \cdot vL_{x,y}\,. \tag{\Al} \eqnlabel{AL}
\]
It follows from \eqnref{AL} that
$(u\ld v)L_{x,y} = u L_{x,y}\ld v L_{x,y}$ and also
$J L_{x,y} = L_{x,y} J$.

The assertion that a permutation $\varphi$ of a loop $Q$ is an automorphism
of $Q$ can be expressed in equivalent ways in terms of the various loop
permutations:
\[
L_x \varphi = \varphi L_{x\varphi}, \qquad
D_x \varphi = \varphi D_{x\varphi}\,.
\]
We shall use these in calculations while referencing \eqnref{AL}.

\section{Preliminaries}
\seclabel{prelims}

In this section, we establish several preliminary results for commutative
A-loops which will be needed later. Some of these generalize rather easily
to arbitrary A-loops, and some of those generalizations can be found in \cite{BP}.
We give brief proofs in the commutative case to make the paper self-contained.

For an automorphism $\varphi$ of a loop $Q$, let
$\Fix(\varphi) = \setof{x\in Q}{x\varphi = x}$. We begin with an easy
observation.

\begin{lemma}
\lemlabel{fix}
Let $Q$ be a loop and let $\varphi\in \Aut(Q)$. Then
\begin{enumerate}[label=\roman*\emph{)}]
\item\quad $\Fix(\varphi)$ is a subloop,
\item\quad If $x\in \Fix(\varphi)$, then $\sbl{x}\leq \Fix(\varphi)$,
\item\quad For each $x\in \Fix(\varphi)$,
\begin{equation}
\eqnlabel{commutes}
L_x \varphi = \varphi L_x  \qquad\text{and}\qquad
D_x \varphi = \varphi D_x \,.
\end{equation}
\end{enumerate}
\end{lemma}

\begin{lemma}
\lemlabel{fix-exchange}
For all $x,y,z$ in a commutative A-loop $Q$,
\[
x\in \Fix(L_{y,z}) \qquad\Leftrightarrow\qquad
y L_x L_z = y L_z L_x \qquad\Leftrightarrow\qquad
z\in \Fix(L_{y,x})\,.
\]
\end{lemma}

\begin{proof}
We have $x L_{y,z} = x$ iff $x L_y L_z = x L_{yz}$
iff $y L_x L_z = y L_z L_x$. Since this last equation is
symmetric in $x$ and $z$, the other equivalence follows.
\end{proof}

For $x$ in a loop $Q$ and $n\in \mathbb{Z}$, we define $x^n = 1 L_x^n$.
Then $x \cdot x^n = 1 L_x^n L_x = 1 L_x^{n+1} = x^{n+1}$ for all
$n\in \mathbb{Z}$. Also, for any $\varphi \in \Aut(Q)$,
$(x^n)\varphi = 1 L_x^n \varphi = 1 \varphi L_{x\varphi}^n = (x\varphi )^n$.

\begin{lemma}[\cite{BP}, Thm 2.6]
\lemlabel{specials1}
In a commutative A-loop, the following identities hold for all $x,y$
and for all $m,n\in \mathbb{Z}$:
\begin{align}
x^n L_{y,x^m} &= x^n  \eqnlabel{fix-x} \\
L_{x^m} L_{x^n} &= L_{x^n} L_{x^m} \eqnlabel{switch-mn} \\
L_{x^n} L_{y,x^m} &= L_{y,x^m} L_{x^n} \eqnlabel{L-commute} \\
D_{x^n} L_{y,x^m} &= L_{y,x^m} D_{x^n} \eqnlabel{D-commute}
\end{align}
\end{lemma}

\begin{proof}
First, we have
$x L_{y,x} = xy \ld (x\cdot yx) = xy\ld (xy\cdot x) = x$, so
that $x\in \Fix(L_{y,x})$. By \eqnref{AL} and Lemma \lemref{fix}(ii),
$x^n \in \Fix(L_{y,x})$ for all $n\in \mathbb{Z}$. Thus by
Lemma \lemref{fix-exchange}, $x\in \Fix(L_{y,x^n})$, and so
$x^m \in \Fix(L_{y,x^n})$ for all $m,n\in \mathbb{Z}$ by
\eqnref{AL} and Lemma \lemref{fix}(ii) again. This establishes
\eqnref{fix-x}, and then \eqnref{switch-mn} follows from
another application of Lemma \lemref{fix-exchange}.
Finally, \eqnref{L-commute} and \eqnref{D-commute} follow from
\eqnref{fix-x} and \eqnref{commutes}.
\end{proof}

A loop is said to be \emph{power-associative} if for each $x$, the
subloop $\sbl{x}$ is a group. Power-associativity is equivalent to
$x^m x^n = x^{m+n}$ for all $x\in Q$ and all $m,n\in \mathbb{Z}$.

\begin{lemma}[\cite{BP}, Thm. 2.4]
\lemlabel{pa}
Every commutative A-loop is power-associative.
\end{lemma}

\begin{proof}
For all $m,k\in \mathbb{Z}$ and for all $x$,
\[
x^m x^{k+1} = x^m (x^k\cdot x) \by{L-commute} x^k (x^m\cdot x) = x^{m+1} x^k\,.
\]
By an easy induction, $x^m x^{k+n} = x^{m+n} x^k$ for all $m,n,k\in \mathbb{Z}$.
Taking $k = -n$, we have the desired result.
\end{proof}

\begin{lemma}
\lemlabel{specials}
In a commutative A-loop, the following identities hold:
\begin{align}
y^n L_{y,x} &= (xy \ld x)^{-n} \qquad\text{ for all } n\in \mathbb{Z}\,, \eqnlabel{fix-y} \\
x y^2 &= (xy)(xy\ld x)\inv\,. \eqnlabel{Omnipresent}
\end{align}
\end{lemma}

\begin{proof}
We compute
\[
y^{-n} L_{y,x} = (y\inv)^n L_{y,x} \by{AL} (y\inv L_{y,x})^n = (xy\ld x)^n\,,
\]
and thus obtain \eqnref{fix-y} upon replacing $n$ with $-n$.
Finally we have
\[
xy \ld xy^2 = y L_{y,x} \by{fix-y} (xy\ld x)\inv\,,
\]
which is equivalent to \eqnref{Omnipresent}.
\end{proof}

A loop is said to have the \emph{automorphic inverse property} (\AIP)
if it has two-sided inverses and satisfies
\begin{equation}
\eqnlabel{AIP} \tag{\AIP}
(xy)\inv = x\inv y\inv
\qquad\text{or equivalently,}\qquad
L_x J = J L_{x\inv}
\end{equation}
for all $x,y$.

\begin{lemma}
\lemlabel{aaip}
Every commutative A-loop has the AIP.
\end{lemma}

\begin{proof}
Using the fact that $L_{x\inv} L_x = L_{x\inv,x}$ is an automorphism, we compute
\begin{alignat*}{2}
y L_x L_{x\inv} J
&\by{switch-mn} y L_{x\inv} L_x J
&&\by{AL} y\inv L_{x\inv} L_x \\
&\ = x\inv \lbrack L_{y\inv} L_y\inv\rbrack \lbrack L_y L_x \rbrack
&&\by{switch-mn} x\inv L_y\inv \lbrack L_{y\inv} L_{y,x}\rbrack L_{xy} \\
&\by{AL} x\inv L_y\inv L_{y,x} L_{y\inv L_{y,x}} L_{xy}
&&\by{fix-y} [(xy)\inv \cdot (xy\ld x)] L_{xy} \\
&\ = x L_{xy}\inv L_{(xy)\inv} L_{xy}
&&\by{switch-mn} x L_{(xy)\inv} \\
&\ = (xy)\inv L_x
&&\ = y L_x J L_x \,.
\end{alignat*}
Thus $L_x L_{x\inv} J = L_x J L_x$, or $L_{x\inv} J = J L_x$.
Replacing $x$ with $x\inv$, we obtain (\AIP).
\end{proof}

\begin{lemma}
\lemlabel{rewrites}
In a commutative A-loop, the following identities hold.
\begin{align}
L_{x,y} &= L_{x\inv, y\inv}  \eqnlabel{aip-linn} \\
L_{x,y} &= L_{x\inv \ld y}\inv L_x L_y  \eqnlabel{rewrite1} \\
L_{x,y} &= L_y L_{x\inv \ld y}\inv L_x \eqnlabel{rewrite2} \\
L_{x\ld y,x} &= L_{(y\ld x)\inv,x}  \eqnlabel{cute} \\
L_{(x\ld y)\inv \ld x}\inv L_{x\ld y} &= L_y\inv L_{y \ld x} \eqnlabel{tripstar}
\end{align}
\end{lemma}

\begin{proof}
First, \eqnref{aip-linn} is an easy consequence of the AIP:
\[
\left( z L_{x,y} \right)\inv \by{AIP} z\inv L_{x\inv,y\inv} \by{AL} \left( z L_{x\inv,y\inv}\right)\inv\,.
\]
For \eqnref{rewrite1}, we compute
\begin{alignat*}{2}
L_{x\inv \ld y}\inv \lbrack L_x L_y \rbrack
&\ = \lbrack L_{x\inv\ld y}\inv L_{x,y}\rbrack L_{yx}
&&\by{AL} L_{x,y} L_{(x\inv\ld y)L_{x,y}}\inv L_{yx} \\
&\by{aip-linn} L_{x,y} L_{(x\inv\ld y)L_{x\inv,y\inv}}\inv L_{yx}
&&\ = L_{x,y} L_{(y\inv x\inv)\inv}\inv L_{yx} \\
&\by{AIP} L_{x,y} L_{yx}\inv L_{yx}
&&\ = L_{x,y}\,.
\end{alignat*}
Next, we have
\[
L_y\inv L_{x,y} \by{L-commute} L_{x,y} L_y\inv
\by{rewrite1} L_{x\inv \ld y}\inv L_x\,,
\]
which gives \eqnref{rewrite2}. For \eqnref{cute}, we
compute
\[
L_{x\ld y,x} = L_{x\ld y,x}\inv L_{x\ld y,x}^2
\by{AL} L_{x\ld y,x}\inv L_{x\ld y,x} L_{(x\ld y)L_{x\ld y,x}, x L_{x\ld y,x}}
= L_{(y\ld x)\inv,x}\,
\]
using \eqnref{fix-y} and \eqnref{fix-x}. Finally,
we apply \eqnref{rewrite1} to both sides of \eqnref{cute} to get
\[
L_{(x\ld y)\inv \ld x}\inv L_{x\ld y} L_x = L_{(y\ld x)\ld x}\inv L_{y\ld x} L_x\,.
\]
Canceling and using $(y\ld x)\ld x = y$, we obtain \eqnref{tripstar}.
\end{proof}

\begin{lemma}
\lemlabel{D-split}
For all $x,y$ in a commutative A-loop,
\begin{align}
D_{x^2} &= D_x J D_x \eqnlabel{D-split} \\
x^2 &= y D_x \cdot y\inv D_x \eqnlabel{D-split2} \\
x &= y\inv D_{x\inv} \cdot y D_{x^2} \eqnlabel{D-split3}\,.
\end{align}
\end{lemma}

\begin{proof}
For all $x,y$,
\begin{alignat*}{3}
y D_{x^2}
&\ = x L_x L_y\inv
&&\ = x L_{x\ld y}\inv \lbrack L_{x\ld y} L_x L_{x\cdot x\ld y}\inv\rbrack
&&\ = x L_{x\ld y}\inv L_{x\ld y,x} \\
&\by{AL} x L_{x\ld y,x} L_{(x\ld y) L_{x\ld y,x}}\inv
&&\by{fix-y} x L_{x\ld y,x} L_{(y\ld x)\inv}\inv
&&\by{fix-x} x L_{(y\ld x)\inv}\inv \\
&\ = (y\ld x)\inv D_x
&&\ = y D_x J D_x\,. &&
\end{alignat*}
This establishes \eqnref{D-split}. Rewrite
\eqnref{D-split} as $J D_x = D_x D_{x^2}$ since
$D_x\inv = D_x$. Applying
this to $y$, we have
$y\inv D_x = y D_x D_{x^2} = x^2 L_{y D_x}\inv$, which is
equivalent to \eqnref{D-split2}. Finally, rewrite
\eqnref{D-split} (applied to $y$) as $x L_{yD_x J}\inv = y D_{x^2}$,
or $x = y D_{x^2} L_{yD_x J}$. Using \eqnref{AIP}, we obtain
\eqnref{D-split3}.
\end{proof}

\section{Commutative A-loops of odd order}
\seclabel{odd-order}

A loop is \emph{uniquely} $2$-\emph{divisible} if the
squaring map $x\mapsto x^2$ is a permutation. In finite,
power-associative loops, being uniquely $2$-divisible is
equivalent to each element having odd order.

The following is well-known and holds in more generality
than we need here.

\begin{lemma}
\lemlabel{odd-2div}
A finite commutative loop $Q$ is
uniquely $2$-divisible if and only if it has odd order.
\end{lemma}

\begin{proof}
If $Q$ is uniquely $2$-divisible, then the
inversion permutation $J$ does not fix any nonidentity elements.
Hence the set of nonidentity elements of $Q$ has even order,
and so $Q$ has odd order.

Now assume $Q$ has odd order, and fix $c\in Q$.
By commutativity, the set $U = \setof{(x,y)}{xy = c, x\neq y}$ has
even order. Since the set $V = \setof{(x,y)}{xy = c}$ has
size $|Q|$, it follows that the set
$U \backslash V = \setof{(x,x)}{x^2 = c}$ has odd order,
and hence is nonempty. Thus the squaring map
$x\mapsto x^2$ is surjective, and hence, by finiteness,
bijective.
\end{proof}

In this section we will study the structure of commutative A-loops
of odd order in detail. To explain our approach, we first need a
useful notion from group theory; \textit{cf.} \cite{Asch2, FKP}.

A \emph{twisted subgroup} of a group $G$ is a subset $T\subset G$
satisfying (i) $1\in T$, (ii)~$a\inv \in T$ for each $a\in T$, and
(iii) $aba\in T$ for each $a,b\in T$. A twisted subgroup $T$ is
uniquely $2$-divisible if the restriction of the squaring map
$x\mapsto x^2$ to $T$ is a permutation.

On a uniquely $2$-divisible twisted subgroup $T$, one can define
a loop operation $\circ$ by $a\circ b = (ab^2 a)^{1/2}$ where
the exponent $1/2$ denotes the unique square root in $T$. The
loop $(T,\circ)$ is then a (left) \emph{Bol loop}, that is,
it satisfies the identity
$x\circ (y\circ (x\circ z)) = (x\circ (y\circ x))\circ z$.
In addition, $(T,\circ)$ satisfies (\AIP); left Bol loops
with (\AIP) are known as left \emph{Bruck loops}.

For some classes of loops, the multiplication groups contain
natural twisted subgroups. Up until now, the only known
example of this is the variety of Bol loops: for a Bol loop
$Q$, the set $L_{Q} = \setof{L_x}{x\in Q}$ of left translations
is a twisted subgroup of $\Mlt(Q)$. In case $Q$ is uniquely
$2$-divisible, there is also a natural left Bruck loop structure
on $L_{Q}$. It turns out that this Bruck loop structure can
be isomorphically transferred to the underlying set $Q$ itself,
so that $Q$ has two loop structures (which may or may not coincide);
its original Bol loop structure and the transferred Bruck loop structure.

There are two things that make all of this particularly useful.
The first is that uniquely $2$-divisible Bruck loops are highly
structured \cite{Glauberman1}. The second is that powers of
elements in the two loop structures coincide. It is thus
possible to prove results about the original Bol loop by
using its associated Bruck loop. This idea was fruitfully exploited
for Moufang loops by Glauberman \cite{Glauberman2}; for the
Bol case, see \cite{FKP}.

\smallskip

We will now apply the same circle of ideas to commutative A-loops.
We will start by identifying a twisted subgroup of the multiplication
group of a commutative A-loop.
For each $x$ in a commutative A-loop $Q$, set
\begin{equation}
\eqnlabel{P}\tag{P}
P_x = L_x L_{x\inv}\inv \by{switch-mn} L_{x\inv}\inv L_x\,.
\end{equation}
and let $P_Q = \setof{P_x}{x\in Q}$. Observe that the set
$P_Q$ trivially satisfies two of the conditions for being
a twisted subgroup: $\id{Q} = P_1\in P_Q$, and for each $x\in Q$,
\[
P_x P_{x\inv} = L_x L_{x\inv}\inv L_{x\inv} L_x\inv = \id{Q}\,,
\]
so that $P_x\inv = P_{x\inv} \in P_Q$.

\begin{lemma}
\lemlabel{twist-tmp}
For all $x,y$ in a commutative A-loop $Q$,
\begin{align}
x\inv P_{xy} &= x y^2  \eqnlabel{twist-tmp2} \\
L_{x\inv} P_{xy} &= P_y L_x \eqnlabel{twist-tmp}
\end{align}
\end{lemma}

\begin{proof}
Applying \eqnref{AIP} to \eqnref{Omnipresent} and rearranging
gives \eqnref{twist-tmp2}. Next, for all $x,y\in Q$,
\begin{alignat*}{3}
L_{x\inv} P_{xy}
&\ = L_{x\inv} L_{(xy)\inv}\inv L_{xy}
&&\by{AIP} L_{x\inv} L_{x\inv y\inv}\inv L_{xy}
&&= L_{y\inv}\inv L_{y\inv,x\inv} L_{xy} \\ 
&\by{aip-linn} L_{y\inv}\inv L_{y,x} L_{xy}
&&\ = L_{y\inv}\inv L_y L_x
&&= P_y L_x\,.
\end{alignat*}
This proves \eqnref{twist-tmp}.
\end{proof}

Note that \eqnref{twist-tmp2} can also be obtained
by applying \eqnref{twist-tmp} to $1\in Q$.

\begin{lemma}
\lemlabel{twisted}
For all $x,y$ in a commutative A-loop $Q$,
\begin{equation}
\eqnlabel{twisted}
P_x P_y P_x = P_{yP_x}\,.
\end{equation}
In particular, $P_Q$ is a twisted subgroup of $\Mlt(Q)$.
\end{lemma}

\begin{proof}
For all $x,y\in Q$,
\begin{align*}
P_x P_y P_x
&&&= P_x P_y L_x L_{x\inv}\inv &&\by{twist-tmp} P_x L_{x\inv} P_{xy} L_{x\inv}\inv \\
&&&\by{P} L_x P_{xy} L_{x\inv}\inv &&= L_x P_{x\inv (x\inv\ld xy)} L_{x\inv}\inv \\
&&&\by{P} L_x P_{x\inv\cdot yP_x} L_{x\inv}\inv &&\by{twist-tmp} P_{yP_x} L_{x\inv} L_{x\inv}\inv \\
&&&= P_{yP_x}\,. &&
\end{align*}
This establishes \eqnref{twisted}, and the rest
follows immediately.
\end{proof}

\begin{lemma}
\lemlabel{power-alt}
For all $x$ in a commutative A-loop $Q$ and for all $n\in \mathbb{Z}$,
\begin{equation}
\eqnlabel{power-alt}
P_x^n = P_{x^n}\,.
\end{equation}
\end{lemma}

\begin{proof}
We have already noted \eqnref{power-alt} for $n=-1$, while it is trivial
for $n = 0,1$. If \eqnref{power-alt} holds some $n$, then
\[
P_x^{n+2} = P_x P_{x^n} P_x \by{twisted} P_{x^n P_x} = P_{x^{n+2}}\,,
\]
the last equality holding by power-associativity (Lemma \lemref{pa}).
The rest follows by induction.
\end{proof}

In calculations, we will frequently use \eqnref{power-alt} without
explicit reference.

Now assume $Q$ is a uniquely $2$-divisible, commutative A-loop. By
\eqnref{power-alt}, the twisted subgroup $P_Q$ is also uniquely $2$-divisible.
Thus there is a natural Bruck loop operation $\circ$ on $P_Q$ given by
\begin{equation}
\eqnlabel{P-bruck}
P_x \circ P_y = ( P_x P_y^2 P_x )^{1/2}
\by{power-alt} ( P_x P_{y^2} P_x )^{1/2}
\by{twisted} ( P_{y^2 P_x} )^{1/2}
\by{power-alt} P_{(y^2 P_x)^{1/2}}\,.
\end{equation}
Thus as with uniquely $2$-divisible Bol loops \cite{FKP} or Moufang loops
\cite{Glauberman2}, we define a new binary operation (for which we
will use the same symbol) on the underlying set $Q$ by
\begin{equation}
\eqnlabel{bruck}\tag{B}
x\circ y = (y^2 P_x)^{1/2} = (x\inv \ld xy^2)^{1/2}\,.
\end{equation}
By \eqnref{P-bruck}, the mapping $x\mapsto P_x$ is a surjective homomorphism
from the magma $(Q,\circ)$ to the loop $(P_Q,\circ)$. In addition, note that
this mapping is injective; indeed, if $P_x = \id{Q}$, then
$x^2 = 1 P_x = 1$ so that $x = 1$. Thus $(Q,\circ)$ is isomorphic to
$(P_Q,\circ)$.
Therefore we have most of the following.

\begin{lemma}
\lemlabel{bruck}
For a uniquely $2$-divisible, commutative A-loop $Q$, $(Q,\circ)$ is
a Bruck loop. Powers in $Q$ coincide with powers in $(Q,\circ)$.
\end{lemma}

\begin{proof}
The remaining assertion about powers follows
easily from \eqnref{bruck}, the power-associativity of $Q$ (Lemma \lemref{pa}),
and an easy induction argument.
\end{proof}

In the finite case, we may now reap the benefits of the known structure
theory of Bruck loops of odd order \cite{Glauberman1} to obtain \emph{Lagrange}
and \emph{Cauchy} theorems. We will implicitly
use Lemma \lemref{odd-2div} in what follows.

\begin{proposition}
\prplabel{lagrange} Let $A\leq B$ be subloops of a finite commutative A-loop
$Q$ of odd order. Then $|A|$ divides $|B|$. In particular, the order of any
element of $Q$ divides $|Q|$.
\end{proposition}

\begin{proof}
The subloops $A$ and $B$ of $Q$ yield subloops $(A,\circ)$ and $(B,\circ)$
of $(Q,\circ)$. The result then follows from (\cite{Glauberman1}, Corollary 4,
p. 395).
\end{proof}

\begin{proposition}
\prplabel{HSC} Let $Q$ be a finite, commutative A-loop of odd order.
If a prime $p$ divides $|Q|$, then $Q$ has an element of order $p$.
\end{proposition}

\begin{proof}
This holds in the corresponding Bruck loop $(Q,\circ)$ \cite{Glauberman1}.
Since powers of an element agree in both $Q$ and $(Q,\circ)$, the result follows.
\end{proof}

\begin{lemma}
\lemlabel{aut-bruck}
Every inner mapping of a uniquely $2$-divisible, commutative A-loop $Q$
acts as an automorphism of $(Q,\circ)$.
\end{lemma}

\begin{proof}
This is obvious from the definition of $\circ$.
\end{proof}

\begin{lemma}
\lemlabel{odd-subloop}
Let $Q$ be a commutative A-loop of odd order. A subloop $K$ of
$(Q,\circ)$ is a subloop of $Q$ if and only if $K \varphi = K$
for each $\varphi \in \Inn(Q)\cap \sbl{L_x : x\in K}$.
\end{lemma}

\begin{proof}
The ``only if'' direction is trivial, so assume the hypothesis
of the converse. Fix $u,v \in K$. Note that $u\inv, v\inv\in K$,
and since powers agree in $(Q,\circ)$ and $Q$,
$v^{1/2}\in K$. Thus $K$ also contains
\[
(u \circ v^{1/2})^2 = v L_u L_{u\inv}\inv
= v L_u^2 L_u\inv L_{u\inv}\inv
= v L_u^2 L_{u\inv,u}\inv\,.
\]
By hypothesis, $K$ then also contains $vL_u^2$.
By induction, $K$ contains $v L_u^{2k}$ for all integers $k$.
Now let $2n+1$ be the order of $u$. Then $L_u^{2n+1} \in \Inn(Q)$,
since $1 L_u^{2n+1} = u^{2n+1} = 1$. Hence
$K$ contains $v L_u^{-2n} L_u^{2n+1} = uv$, and
also $v L_u^{2(-n-1)} L_u^{2n+1}  = u \ld v$. Thus $K$ is closed
under multiplication and left division in $Q$ and is therefore
a subloop of $Q$.
\end{proof}

At a particular point in the proof of Theorem \thmref{ft} below, we
will show that the Bruck loop associated to a certain
commutative A-loop is commutative.
In order to proceed, we will then need the corollary to the following
technical lemma.

\begin{lemma}
\lemlabel{P-commute}
Let $Q$ be a commutative A-loop and assume that the identity
\begin{equation}
\eqnlabel{P-commute}
y^2 P_x = x^2 P_y
\end{equation}
holds for all $x,y\in Q$. Then for all $x,y\in Q$,
\begin{equation}
\eqnlabel{P-iso}
y^2 P_x = x^2 y^2\,.
\end{equation}
\end{lemma}

\begin{corollary}
\corlabel{P-commute}
Let $Q$ be a uniquely $2$-divisible, commutative A-loop. Then
$(Q,\circ)$ is commutative if and only if $(Q,\circ)$ is
isomorphic to $Q$.
\end{corollary}

Indeed, in the uniquely $2$-divisible case, \eqnref{P-commute}
asserts that $(Q,\circ)$ is commutative, and
\eqnref{P-iso}
says that $(x\circ y)^2 = x^2 y^2$, that is, the
squaring map $x\mapsto x^2$ is an isomorphism from
$(Q,\circ)$ to $Q$.

\begin{proof}[Proof of Lemma \lemref{P-commute}]
First we establish
\begin{equation}
\eqnlabel{P-comm-tmp3}
(xy^2) P_x = x P_{xy}
\end{equation}
for all $x,y\in Q$. Indeed, we have
\begin{alignat*}{4}
(xy^2) P_x
&\by{switch-mn} y^2 P_x L_x
&&\by{P-commute} x^2 P_y L_x
&&\by{P} 1 P_x P_y P_x L_{x\inv}
&&\by{twisted} 1 P_{y P_x} L_{x\inv} \\
&\ = x\inv  \left( y P_x \right)^2
&&\by{twist-tmp2} x P_{x\inv \cdot yP_x}
&&\by{P} x P_{xy}\,. &&
\end{alignat*}

Next, we will also require
\begin{equation}
\eqnlabel{P-comm-tmp4}
x\inv P_{y^2} = y^2 P_{x\inv} P_y L_x
\end{equation}
for all $x,y\in Q$. For this, we compute
\begin{alignat*}{2}
x\inv P_{y^2}
&= x\inv P_{x\cdot x\ld y^2}
&&\by{twist-tmp2} x (x\ld y^2)^2 \\
&= (x\ld y^2)^2 P_{y\inv} P_y L_x
&&\by{P-commute} y^{-2} P_{x\ld y^2} P_y L_x \\
&= \left( (x\ld y^2) \cdot (x\ld y^2)\inv \ld y^{-2} \right) P_y L_x
&&\by{AIP} \left( (x\ld y^2) \cdot (x\inv \ld y^{-2}) \ld y^{-2} \right) P_y L_x \\
&= \left( (x\ld y^2) x\inv \right) P_y L_x
&&\ = y^2 P_{x\inv} P_y L_x\,.
\end{alignat*}

Now, we compute
\begin{alignat*}{2}
y^2 P_x P_y L_x
&\by{P-commute} x^2 P_y^2 L_x
&&\by{P} 1 P_x P_{y^2} P_x L_{x\inv} \\
&\by{twisted} 1 P_{y^2 P_x} L_{x\inv}
&&\ = x\inv \left( y^2 P_x \right)^2 \\
&\by{twist-tmp2} x P_{x\inv\cdot y^2 P_x}
&&\ = x P_{xy^2} \\
&\by{AIP} x P_{x\inv y\inv} P_{xy} P_{xy^2}
&&\by{twist-tmp2} \left(x\inv y^{-2}\right) P_{xy} P_{xy^2} \\
&\by{AIP} (xy^2)\inv P_{xy^2\cdot (xy^2 \ld xy)} P_{xy^2}
&&\by{twist-tmp2} \left( xy^2 \cdot (xy^2 \ld xy)^2\right) P_{xy^2}\\
&\by{P-comm-tmp3} (xy^2) P_{xy^2 \cdot (xy^2 \ld xy)}
&&\ = (xy^2) P_{xy} \\
&\by{twist-tmp2} x\inv P_{xy}^2
&&\ = x\inv P_{(xy)^2} \\
&\by{P-comm-tmp4} (xy)^2 P_{x\inv} P_{xy} L_x\,. &&
\end{alignat*}
Canceling $L_x$, we have
\[
y^2 P_x P_y  = (xy)^2 P_{x\inv} P_{xy} = 1 P_{xy} P_{x\inv} P_{xy}
= 1 P_{x\inv P_{xy}} \by{twist-tmp2} 1 P_{xy^2} = (xy^2)^2\,.
\]
Thus
\[
y^2 P_x = (xy^2)^2 P_{y\inv} \by{P-commute} y^{-2} P_{y^2 x}
\by{twist-tmp2} y^2 x^2\,,
\]
which is \eqnref{P-iso}.
\end{proof}

We now turn to the main result of this section

\begin{theorem}[Odd Order Theorem]
\thmlabel{ft}
Every commutative A-loop of odd order is solvable.
\end{theorem}

\begin{proof}
Let $Q$ be a minimal counterexample. Since normal subloops
and quotients of commutative A-loops of odd order also have
odd order, it follows that $Q$ must be simple.
Let $N$ denote the derived subloop of $(Q,\circ)$, that is, the
smallest normal subloop of $(Q,\circ)$ such that $(Q/N,\circ)$
is an abelian group. Finite Bruck loops of odd order are solvable
(\cite{Glauberman2}, Thm. 14(b)), and so $N$ is a proper subloop.
Clearly $N$ is fixed by every automorphism of $(Q,\circ)$.
By Lemma \lemref{aut-bruck}, $N$ is fixed by every element of
$\Inn(Q)$. Thus by Lemma \lemref{odd-subloop}, $N$ is a subloop
of $Q$ itself. Since $N$ is invariant under $\Inn(Q)$, $N$ is
normal in $Q$. But $Q$ is simple, and so $N = \{1\}$. Therefore
$(Q,\circ)$ is an abelian group. By Corollary \corref{P-commute},
$(Q,\circ)$ is isomorphic to $Q$. Thus $Q$ is an abelian group,
which contradicts the assumption that $Q$ is not solvable.
\end{proof}

\section{Squares and an Associated Loop}
\seclabel{squares}

In an abelian group, or even a commutative Moufang loop, the product
of two squares is trivially a square, for in such loops the identity
$x^2 y^2 = (xy)^2$ holds. This identity does not hold in commutative
A-loops. For example, there is a nonassociative, commutative A-loop
of order $15$ \cite{Drapal} in which the identity fails.
Nevertheless, the more fundamental assertion about the product of two
squares holds, as we are going to show.

Motivated by Theorem \thmref{SquareRoot} below,
we introduce a new binary operation in commutative A-loops:
\[
\eqnlabel{newop}
    x \nlp y = \left( xy\ld x \cdot yx\ld y \right)\inv = y L_{y,x} \cdot x L_{x,y}\,, \tag{$\nlp$}
\]
where the second equality follows from \eqnref{fix-y} and \eqnref{AIP}.

\begin{theorem}
\thmlabel{SquareRoot}
For all $x,y$ in a commutative A-loop,
\[
    x^2 y^2 = \left( x \nlp y \right)^2\,.
\]
\end{theorem}

To establish the theorem, we require a couple of lemmas.

\begin{lemma}
\lemlabel{perms}
For all $x,y$ in a commutative A-loop $Q$,
\begin{equation}
\eqnlabel{perms}
x \nlp y = x^2 \cdot x\ld (xy \ld x)\inv \,.
\end{equation}
\end{lemma}

\begin{proof}
First, we have
\begin{equation}
\eqnlabel{perms1}
x L_{x,y}
= (x^2 y)L_{yx}\inv
= y L_x\inv L_x L_{x^2} L_{yx}\inv
\by{switch-mn} y L_x\inv L_{x^2} L_x L_{yx}\inv
= y L_x\inv L_{x^2} L_y\inv L_{y,x}\,.
\end{equation}
Thus,
\begin{alignat*}{2}
x \nlp y
&\ = y L_{y,x} \cdot x L_{x,y}
&&\by{perms1} y L_{y,x}\cdot y L_x\inv L_{x^2} L_y\inv L_{y,x} \\
&\by{AL} \lbrack y \cdot y L_x\inv L_{x^2} L_y\inv \rbrack L_{y,x}
&&\ = y L_x\inv L_{x^2} L_{y,x} \\
&\by{commutes} y L_{y,x} L_x\inv L_{x^2}
&&\by{fix-y} (xy\ld x)\inv L_x\inv L_{x^2} \\
&\ = x^2 \cdot x\ld (xy \ld x)\inv\,, &&
\end{alignat*}
which gives \eqnref{perms}.
\end{proof}

\begin{lemma}
\lemlabel{newfix}
For all $x,y$ in a commutative A-loop,
\begin{equation}
\eqnlabel{strangecomm}
x\inv \ld (xy \ld x) = y\ld (yx \ld y)\inv\,.
\end{equation}
\end{lemma}

\begin{proof}
We compute
\begin{alignat*}{2}
(y\ld (yx\ld y)\inv) L_{x\inv} L_{xy}
&\ = (yx\ld y)\inv L_{x \ld x y}\inv L_{x\inv} L_{xy}
&&\by{rewrite1} \left(yx\ld y\right)\inv L_{x\inv,xy} \\
&\by{AL} \left( (xy)L_{x\inv,xy} \ld y L_{x\inv,xy} \right)\inv
&&\by{fix-x} \left( xy \ld y L_{x\inv,xy} \right)\inv \\
&\by{aip-linn} \left( xy \ld y L_{x,(xy)\inv} \right)\inv
&&\ = \left( xy \ld (x (xy)\inv)\inv \right)\inv \\
&\by{AIP} \left( xy \ld (x\inv \cdot xy) \right)\inv
&&\ = x\,.
\end{alignat*}
Thus $y\ld (yx\ld y)\inv = x L_{xy}\inv L_{x\inv}\inv = x\inv \ld (xy \ld x)$, as claimed.
\end{proof}

Now we turn to the main result of this section.

\begin{proof}[Proof of Theorem \thmref{SquareRoot}]
Set $z = x \nlp y$. Then
\begin{alignat*}{2}
x^2 D_z
&\ = z L_{x^2}\inv
&&\by{perms} (x^2 \cdot x\ld (xy\ld x)\inv))L_{x^2}\inv \\
&\ = x \ld (xy\ld x)\inv
&&\by{AIP} (x\inv \ld (xy\ld x))J  \\
&\by{strangecomm} (y \ld (yx\ld y)\inv) J
&&\ = ( y^2 \cdot y\ld (yx\ld y)\inv ) L_{y^2}\inv J  \\
&\by{perms} z L_{y^2}\inv J
&&\ = y^2 D_z J\,.
\end{alignat*}
Thus $x^2 = x^2 D_z^2 = y^2 D_z J D_z  \by{D-split} y^2 D_{z^2} = z^2 L_{y^2}\inv$, and so
$x^2 y^2 = z^2$, as claimed.
\end{proof}

As the notation suggests, we will now consider $(Q,\nlp)$ as being a
new magma constructed on a commutative A-loop $Q$.
We introduce notation for the corresponding left translation map:
\[
y S_x = x \nlp y  \tag{S}
\]
for all $x,y$. Note that
\begin{equation}
\eqnlabel{S-perm}
S_x = L_x D_x J L_x\inv L_{x^2}
\end{equation}
by Lemma \lemref{perms}.

\begin{proposition}
\prplabel{newop-comm-loop}
Let $Q$ be a commutative A-loop and let $\nlp$ be defined by \eqnref{newop}.
Then $(Q,\nlp)$ is a power-associative, commutative loop with the same
neutral element as $Q$. Powers in $(Q,\nlp)$ coincide with powers in $Q$.
\end{proposition}

\begin{proof}
Commutativity is clear from the definition as is the fact that $(Q,\nlp)$
has the same neutral element as $Q$. By \eqnref{S-perm},
each $S_x$ is a permutation of $Q$.
Hence $(Q,\nlp)$ is a loop. Finally, power-associativity of $(Q,\nlp)$ and the
coinciding of powers follow from the power-associativity of $Q$
(Lemma \lemref{pa}).
\end{proof}

For later use, we note the following.

\begin{lemma}
\lemlabel{S-commute}
For all $x,y$ in a commutative A-loop $Q$ and all $m,n\in \mathbb{Z}$,
\begin{equation}
\eqnlabel{S-commute}
S_{x^n} L_{y,x^m} = L_{y,x^m} S_{x^n}\,.
\end{equation}
\end{lemma}

\begin{proof}
This follows immediately from \eqnref{S-perm}, \eqnref{L-commute}, \eqnref{D-commute}
and (\AIP).
\end{proof}

We conclude this section by noting that
for uniquely $2$-divisible, commutative A-loops,
the loop operation $\nlp$ gives nothing new.

\begin{lemma}
\lemlabel{2div-iso}
If $Q$ is a uniquely $2$-divisible, commutative A-loop, then $(Q,\nlp)$ is
isomorphic to $Q$.
\end{lemma}

\begin{proof}
Indeed, the conclusion of Theorem \thmref{SquareRoot} shows that the squaring map is an
isomorphism from $(Q,\nlp)$ to $Q$.
\end{proof}

We will return to the associated loop operation $(Q,\nlp)$ in
\S\secref{Exp2} when we consider commutative A-loops of
exponent $2$.

\section{The Decomposition Theorem}
\seclabel{decomposition}

Our main goal in this section is the following.

\begin{theorem}[Decomposition for Finite Commutative A-loops]
\thmlabel{Decomposition}
If $Q$ is a finite commutative A-loop, then $Q = K(Q)\times H(Q)$,
where $K(Q) = \setof{x\in Q}{|x|\text{ is odd }}$ and
$H(Q) = \setof{x\in Q}{x^{2^n} = 1\text{ for some }n\in \mathbb{Z}}$.
\end{theorem}

In addition, $K(Q)$ has odd order (Theorem \thmref{IteratedSquares}(v) below),
and we will show later that $H(Q)$ has order a power of $2$
(Theorem \thmref{p-loop}).

\begin{proposition}
\prplabel{Squares}
In a commutative A-loop $Q$, the set $K_1(Q) = \setof{x^2}{x\in Q}$ is a
normal subloop of $Q$.
\end{proposition}

\begin{proof}
The set $K_1$ is closed under multiplication by Theorem \thmref{SquareRoot}.
By Proposition \prpref{newop-comm-loop}, given $x,z\in Q$, there
exists a unique $y\in Q$ such that $x \nlp y = z$, and so $x^2 y^2=z^2$
by Theorem \thmref{SquareRoot} once more. Thus $K_1$ is a subloop of $Q$.
The normality of $K_1$ follows from the fact that all inner mappings of $Q$ are
automorphisms of $Q$ and hence preserve squares.
\end{proof}

\begin{theorem}
\thmlabel{IteratedSquares}
Let $Q$ be a commutative A-loop. For $n\geq 1$, define
\begin{align*}
    K_n(Q) &= \setof{x^{2^n}}{x\in Q},\\
    K(Q) &= \bigcap_{n\geq 1} K_n(Q).
\end{align*}
Then:
\begin{enumerate}[label=\roman*)]
\item\qquad $K_{n+1}(Q) = \setof{x^2}{x\in K_n(Q)}$ for every $n\geq 0$.
\item\qquad $K_{n+1}(Q)\subseteq K_n(Q)$ for every $n\geq 0$.
\item\qquad $K_n(Q)\unlhd Q$ for every $n\geq 0$.
\item\qquad $K(Q)\unlhd Q$.
\item\qquad If $Q$ is finite, then $K(Q) = \setof{x\in Q}{|x|\text{ is odd }}$ and $|K(Q)|$ is odd.
\end{enumerate}
\end{theorem}

\begin{proof}
If $x\in K_n(Q)$ then $x=y^{2^n}$ for some $y\in Q$ and $x^2 = y^{2^{n+1}}\in
K_{n+1}(Q)$. Conversely, if $x\in K_{n+1}(Q)$ then $x=z^{2^{n+1}} =
(z^{2^n})^2$ for some $z\in Q$ and $z^{2^n}\in K_n(Q)$. This proves (i) and
(ii).

By Proposition \prpref{Squares}, $K_1(Q)\leq Q$. Assume that $K_n(Q)\leq Q$. By
(i), Proposition \prpref{Squares} applied to $K_n(Q)$ yields $K_{n+1}(Q)\leq K_n(Q)\leq Q$. The normality of $K_n(Q)$ in the A-loop $Q$ follows for free.
This proves (iii) and (iv).

For (v), assume that $Q$ is finite. Then there is $n$ such
that $K_{n+1}(Q) = K_n(Q) = K(Q) = \setof{x^2}{x\in K(Q)}$, by (i). The mapping
$x\mapsto x^2$ is a bijection of $K(Q)$ fixing $1\in K(Q)$, so $K(Q)$ contains
no elements of order $2$ and hence no elements of even order. Conversely, pick
$x\in Q$ of odd order, say $|x|=2m+1$. The equality $x=x^{2m+2}=(x^{m+1})^2$
then implies $x\in K_1(Q)$, so that $x^{m+1}\in K_1(Q)$ by (iii).
Thus $x\in K_2(Q)$ by (i), and so on, proving $x\in K(Q)$.
The remaining assertion follows from Lemma \lemref{odd-2div}.
\end{proof}

\begin{lemma}
\lemlabel{Left}
For every $x,y$ in a commutative A-loop $Q$,
\begin{equation}
\eqnlabel{weird-square}
(x\ld (y\ld x))^2 \ld (y\inv (y\ld x))^2 = (x\ld y)^{-2}
\end{equation}
\end{lemma}

\begin{proof}
With $y$ replaced by $x\ld y$, \eqnref{Omnipresent} yields
\begin{equation}
\eqnlabel{Left1}
    x (x\ld y)^2 = y(y\ld x)\inv \,.
\end{equation}
Replacing $y$ with $y\ld x$ and using $(y\ld x)\ld x = y$ gives
\begin{equation}
\eqnlabel{iv}
x (x\ld (y\ld x))^2 = y\inv (y\ld x)\,.
\end{equation}
 Applying $J$ and using \eqnref{AIP} gives
\begin{equation}
\eqnlabel{Left2}
x\inv (x\ld (y\ld x))^{-2} = y (y\ld x)\inv \,.
\end{equation}
Putting \eqnref{Left1} and \eqnref{Left2} together, we have
\[
    (x\ld y)^2 (x\ld(y\ld x))^{-2}
    =  x D_{y(y\ld x)\inv}\cdot x\inv D_{y(y\ld x)\inv}
    \by{D-split2} (y(y\ld x)\inv)^2 \,.
\]
Applying $J$ to both sides and using \eqnref{AIP}, we have
$(x\ld y)^{-2} (x\ld(y\ld x))^2= (y\inv (y\ld x))^2$, and
this is clearly equivalent to \eqnref{weird-square}.
\end{proof}

\begin{proposition}\label{Pr:Key}
Let $Q$ be a commutative A-loop, and let $x\in Q$ satisfy $x^{2^n}=1$. Then
$(xy)^{2^n}=y^{2^n}$ for every $y\in Q$.
\end{proposition}

\begin{proof}
We proceed by induction on $n$. The claim is clearly true when $n=0$. Let $n\geq 0$, assume that the claim holds for $n$, and let $x\in Q$ satisfy
$x^{2^{n+1}}=1$. Then the induction assumption yields
\begin{equation}
\eqnlabel{Ind1}
    (x^2y)^{2^n} = y^{2^n} = (x^2(x^2\ld y))^{2^n} = (x^2\ld y)^{2^n}
\end{equation}
for every $y\in Q$. We may apply any automorphism $\varphi$ to
\eqnref{Ind1}, and then set $z = y \varphi$ to obtain
$((x\varphi)^2 z)^{2^n} = z^{2^n} = ((x\varphi)^2 \ld z)^{2^n}$
for all $z\in Q$. In particular, we choose $\varphi = J L_{x,x\ld y}$
(by \eqnref{AL} and \eqnref{AIP}).
Then $x J L_{x,x\ld y} = y\ld (x\ld y)$ by \eqnref{fix-y}
(or direct calculation). Hence
\begin{equation}
\eqnlabel{Ind2}
    (z(y\ld(x\ld y))^2)^{2^n} = z^{2^n} = ((y\ld(x\ld y))^2\ld z)^{2^n}
\end{equation}
for every $y, z\in Q$. Thus
\begin{align*}
    y^{2^{n+1}}
    &\by{Ind2} [y(y\ld(x\ld y))^2]^{2^{n+1}}
    \by{iv} [x\inv(x\ld y)]^{2^{n+1}}
    = [(x\inv(x\ld y))^2]^{2^n}\\
    &\by{Ind2} [(y\ld(x\ld y))^2 \ld (x\inv(x\ld y))^2]^{2^n}
    \by{weird-square} (y\ld x)^{-2^{n+1}}\,.
\end{align*}
Then
\begin{align*}
    (y\inv)^{-2^{n+1}}
    & = y^{2^{n+1}} \by{fix-x} y^{2^{n+1}}  L_{y,y\inv}
    = (y\ld x)^{-2^{n+1}} L_{y,y\inv}\\
    &\by{AL} ((y\ld x) L_{y,y\inv})^{-2^{n+1}}
    = (y\inv x)^{-2^{n+1}}\,.
\end{align*}
Taking inverses and replacing $y$ with $y\inv$, we obtain
$y^{2^{n+1}} = (xy)^{2^{n+1}}$, which completes the proof.
\end{proof}

\begin{theorem}
\thmlabel{Power2Exp}
Let $Q$ be a commutative A-loop. For $n\geq 0$, let
\begin{align*}
    H_n(Q) &= \setof{x\in Q}{x^{2^n}=1},\\
    H(Q) &= \bigcup_{n\geq 0} H_n(Q).
\end{align*}
Then:
\begin{enumerate}[label=\roman*\emph{)}]
\item\quad $H_{n+1}(Q) = \setof{x\in Q}{x^2\in H_n(Q)}$ for every $n\geq 0$.
\item\quad $H_{n+1}(Q) \supseteq H_n(Q)$ for every $n\geq 0$.
\item\quad $H_n(Q)\unlhd Q$ for every $n\geq 0$.
\item\quad $H(Q)\unlhd Q$.
\end{enumerate}
\end{theorem}

\begin{proof}
Parts (i) and (ii) are obvious. For (iii) and (iv), it suffices to show that
$H_n(Q)\leq Q$ for every $n\geq 0$ and $H(Q)\leq Q$. Let $x\in H_n(Q)$, $y\in
H_m(Q)$ and let $k=\max\{n,m\}$. Then Proposition \ref{Pr:Key} yields
$(xy)^{2^k} = x^{2^k} = 1$ and $(x\ld y)^{2^k} = (x\cdot x\ld y)^{2^k} =
y^{2^k} = 1$.
\end{proof}

Finally, we turn to the proof of the main result of this section.

\begin{proof}[Proof of Theorem \thmref{Decomposition}]
By Theorems \thmref{IteratedSquares} and \thmref{Power2Exp},
$K$ and $H$ are normal subloops of $Q$.
Clearly $K\cap H=1$, and $K H=Q$ is proved in
the same way as for groups (since the argument takes place
in cyclic subgroups, by power-associativity). Then
$Q = K\times H$ follows.
\end{proof}

\section{Commutative A-loops of exponent $2$}
\seclabel{Exp2}

We now turn to commutative A-loops of exponent $2$.
The following result shows why this special case
is of particular importance.

\begin{proposition}
\prplabel{simple}
A finite simple commutative A-loop is either a cyclic group
of order $p$ for some odd prime $p$, or it has exponent $2$.
\end{proposition}

\begin{proof}
Let $Q$ be a finite simple commutative A-loop.
By the Decomposition Theorem \thmref{Decomposition},
$Q = K(Q)\times H(Q)$. Since $Q$ is simple,
$Q = K(Q)$ or $Q = H(Q)$. In the former case, $Q$
is solvable by Theorems \thmref{IteratedSquares}(v) and
\thmref{ft}. Thus $Q$ is both simple and solvable, and
hence is a cyclic group of odd prime order.
Now assume $Q = H(Q)$, that is, every element of $Q$
has order a power of $2$. The subloop
$K_1(Q) = \setof{x^2}{x\in Q}$
is normal (Proposition \prpref{Squares}), and so
either $K_1(Q) = Q$ or $K_1(Q) = \sbl{1}$. In the
former case, the squaring map is a bijection by finiteness,
but then $Q$ has odd order by Lemma \lemref{odd-2div},
a contradiction. Thus for every $x\in Q$, $x^2 = 1$,
that is, $Q$ has exponent $2$.
\end{proof}

Our goal in this section is to establish the following.

\begin{theorem}
\thmlabel{abelian-group}
Let $Q$ be a commutative A-loop of exponent $2$.
Then $(Q,\nlp)$ is an elementary abelian $2$-group.
\end{theorem}

\begin{corollary}
\corlabel{exp2-order}
If $Q$ is a finite, commutative A-loop of exponent $2$,
then $|Q|$ is a power of $2$.
\end{corollary}

The proof of Theorem \thmref{abelian-group} will require
some technical lemmas. \emph{Throughout the rest of this section},
let $Q$ be a commutative A-loop of exponent $2$.
The operation $\nlp$ and the corresponding
translations $S_x$ simplify accordingly:
\begin{align*}
x \nlp y &= x \ld (xy \ld x)  \\
S_x &= L_x D_x L_x\inv
\end{align*}
Thus
$S_x^2 = L_x D_x L_x\inv L_x D_x L_x\inv =
L_x D_x^2 L_x\inv = \id{Q}$. This establishes
the following.

\begin{lemma}
\lemlabel{steiner}
For all $x,y\in Q$, $x \nlp (x \nlp y) = y$, that is,
$S_x^2 = \id{Q}$.
\end{lemma}

\begin{lemma}
\lemlabel{move-em}
For all $x\in Q$,
\begin{equation}
\eqnlabel{move-em}
S_x = L_x D_x L_x\inv = L_x\inv D_x L_x\,.
\end{equation}
\end{lemma}

\begin{proof}
The first equality has already been established. Since $Q$ has
exponent $2$, $D_x = D_{x L_x^2}$ for each $x$. Now $L_x^2 = L_{x,x}\in \Inn(Q)$,
and so we have
$L_x^2 D_x = L_x^2 D_{x L_x^2} \by{AL} D_x L_x^2$. Applying $L_x\inv$ on the
left and on the right, we obtain the desired result.
\end{proof}

\begin{lemma}
\lemlabel{box}
For all $x,y,z\in Q$,
\begin{equation}
\eqnlabel{box}
y L_{z\ld (x\cdot zy),z} S_{zy} = z L_y L_x\inv D_y L_x\,.
\end{equation}
\end{lemma}

\begin{proof}
First, we compute
\begin{alignat*}{2}
y L_{x,z} S_{zy} L_{zx\ld zy}\inv L_{zx}
&\ =  y L_{x,z} S_{zy} L_{zy}\inv \lbrack L_{zy} L_{zx\ld zy}\inv L_{zx}\rbrack
&&\by{rewrite2} y L_{x,z} S_{zy} \lbrack L_{zy}\inv L_{zx,zy}\rbrack \\
&\by{L-commute} y L_{x,z} \lbrack S_{zy} L_{zx,zy}\rbrack L_{zy}\inv
&&\by{S-commute} y \lbrack L_{x,z} L_{zx,zy}\rbrack S_{zy} L_{zy}\inv \\
&\ = \lbrack y L_x\rbrack L_z \lbrack L_{zx}\inv L_{zx}\rbrack  L_{zy} L_{zy\cdot zx}\inv S_{zy} L_{zy}\inv
&&\ = x L_y L_z L_{zy} L_{zy\cdot zx}\inv \lbrack S_{zy} L_{zy}\inv \rbrack \\
&\by{move-em} x L_y L_z L_{zy} L_{zy\cdot zx}\inv L_{zy}\inv D_{zy}\,. &&
\end{alignat*}
Now since $Q$ has exponent $2$, $1 L_y L_z L_{zy} = 1$, and so
$L_y L_z L_{zy}\in \Inn(Q)$. Also,
$zx\cdot zy = (y\ld x)L_y L_z L_{zy}$.
Thus we may apply \eqnref{AL} to get
\begin{alignat*}{2}
y L_{x,z} S_{zy} L_{zx\ld zy}\inv L_{zx}
&= x L_{y\ld x}\inv L_y L_z \lbrack L_{zy} L_{zy}\inv\rbrack D_{zy}
&&= \lbrack x L_{y\ld x}\inv \rbrack L_y L_z D_{zy} \\
&= \lbrack y D_x^2 L_y L_z\rbrack D_{zy}
&&= z D_{zy} \\
&= y\,. &&
\end{alignat*}
where we have used $y^2 = 1$ in the penultimate step. Hence
\[
y L_{x,z} S_{zy}   = y L_{zx}\inv L_{zx\ld zy}
\by{tripstar} y L_{(zy\ld zx)\ld zy}\inv L_{zy\ld zx} \,.
\]
Replacing $x$ with $x L_{zy} L_z\inv = z\ld (x\cdot zy)$, we
obtain
\[
y L_{z\ld (x\cdot zy),z} S_{zy} = y L_{x\ld zy}\inv L_x = z L_y L_x\inv D_y L_x\,.
\]
This establishes \eqnref{box}.
\end{proof}

\begin{lemma}
\lemlabel{star}
For all $u,v,w\in Q$,
\begin{equation}
\eqnlabel{star}
u L_{v\ld (w\cdot uv),v} = u L_v L_w\inv D_v L_w \,.
\end{equation}
\end{lemma}

\begin{proof}
We compute
\begin{alignat*}{2}
u L_{v\ld (w\cdot uv),v}
&\ = \lbrack u L_{v\ld (w\cdot uv)}\rbrack L_v L_{w\cdot uv}\inv
&&\ = w \lbrack L_{uv} L_v\inv L_u\rbrack L_v L_{w\cdot uv}\inv \\
&\by{rewrite2} w L_{u,uv} L_v L_{w\cdot uv}\inv
&&\ = w L_{v\ld uv,uv} L_v L_{w\cdot uv}\inv \\
&\by{cute} w \lbrack L_{uv\ld v,uv} L_v \rbrack L_{w\cdot uv}\inv
&&\ = \lbrack w L_{uv\ld v}\rbrack L_{uv} L_{w\cdot uv}\inv \\
&\ = (uv\ld v)L_w L_{uv}  L_{w\cdot uv}\inv
&&\ = v L_{uv}\inv L_{w,uv} \\
&\by{rewrite2} v L_{w\ld uv}\inv L_w
&&\ = u L_v L_w\inv D_v L_w\,,
\end{alignat*}
which establishes \eqnref{star}.
\end{proof}

\begin{lemma}
\lemlabel{151}
For all $u,v,w\in Q$,
\begin{equation}
\eqnlabel{151}
u L_{v\ld w}\inv L_v L_{vw,u} = wu\,.
\end{equation}
\end{lemma}

\begin{proof}
We compute
\begin{alignat*}{2}
u \lbrack L_{v\ld w}\inv L_v\rbrack L_{vw,u}
&\by{rewrite1} u L_{v,w} L_w\inv L_{vw,u}
&&\by{L-commute}  u L_w\inv \lbrack L_{v,w} L_{vw,u}\rbrack \\
&\ = \lbrack u L_w\inv L_v\rbrack L_w L_u L_{vw\cdot u}\inv
&&\ = v L_{w\ld u} L_w L_u L_{vw\cdot u}\inv \\
&\ = v \lbrack L_{w\ld u} L_{w,u}\rbrack L_{wu} L_{vw\cdot u}\inv
&&\by{rewrite1} v L_w L_u L_{wu} L_{vw\cdot u}\inv \\
&\ = ((u\cdot vw)\cdot wu) L_{vw\cdot u}\inv
&&\ = wu\,,
\end{alignat*}
which establishes \eqnref{151}.
\end{proof}

\begin{lemma}
\lemlabel{209}
For all $u,v,w\in Q$,
\begin{equation}
\eqnlabel{209}
v L_{w,u} S_{uv} = v L_{w,u} L_u\inv L_v\,.
\end{equation}
\end{lemma}

\begin{proof}
We begin with
\[
v L_{u\ld (w\cdot uv),u} S_{uv}
\by{box} u L_v L_w\inv D_v L_w
\by{star} u L_{v\ld (w\cdot vu),v}\,.
\]
Replacing $w$ with $w L_{uv}\inv L_u$, we have
\begin{alignat*}{2}
v L_{w,u} S_{uv}
&\ = u L_{v\ld uw,v}
&&\by{cute} u L_{uw\ld v,v} \\
&\by{rewrite1} u L_{(uw\ld v)\ld v}\inv L_{uw\ld v} L_v
&&\ = u L_{uw}\inv L_{uw\ld v} L_v \\
&\ = v L_{uw}\inv L_{uw\ld u} L_v
&&\by{tripstar} v L_{(u\ld uw)\ld u}\inv L_{u\ld uw} L_v \\
&\ = v L_{w\ld u}\inv L_w L_v
&&\by{rewrite1} v L_{w,u} L_u\inv L_v\,.
\end{alignat*}
This establishes \eqnref{209}.
\end{proof}

\begin{lemma}
\lemlabel{almost}
For all $x,y\in Q$,
\begin{equation}
\eqnlabel{almost}
L_x\inv D_y L_x = L_y\inv D_x L_y D_{xy}\,.
\end{equation}
\end{lemma}

\begin{proof}
We have
\begin{alignat*}{2}
z L_y L_x\inv D_y L_x
&\by{box} y L_{z\ld (x\cdot zy)} S_{zy}
&&\by{star} y \lbrack L_z L_x\inv\rbrack D_z L_x S_{zy} \\
&\ = y L_{z\ld x}\inv \lbrack L_{z\ld x,z} D_z\rbrack L_x S_{zy}
&&\by{D-commute} y L_{z\ld x}\inv D_z \lbrack L_{z\ld x,z} L_x\rbrack S_{zy} \\
&\ = y L_{z\ld x}\inv D_z L_{z\ld x} L_z S_{zy}\,. &&
\end{alignat*}
Now set $u = y L_{z\ld x}\inv D_z L_{z\ld x} =
z L_{(z\ld x)\ld y}\inv L_{z\ld x}$, and observe that
\begin{equation}
\eqnlabel{al-tmp}
u L_{(z\ld x)y,z} \by{151} yz\,.
\end{equation}
Thus using the commutativity of $\nlp$, we compute
\begin{alignat*}{3}
z L_y L_x\inv D_y L_x
&\ = (zu) S_{zy}
&&\ = (zy) S_{zu}
&&\by{al-tmp} u L_{(z\ld x)y,z} S_{zu} \\
&\by{209} u L_{(z\ld x)y,z} L_z\inv L_u
&&\by{al-tmp} (yz) L_z\inv L_u
&&\ = y L_u \\
&\ = u L_y
&&\ = z L_{(z\ld x)\ld y}\inv L_{z\ld x} L_y
&&\by{rewrite1} z L_{z\ld x,y} \\
&\ = z L_{z\ld x} L_y L_{(z\ld x)y}\inv
&&\ = (yx) L_{(z\ld x)y}\inv
&&\ = z D_x L_y D_{xy}\,.
\end{alignat*}
Thus $L_y L_x\inv D_y L_x = D_x L_y D_{xy}$.
Multiplying on the left by $L_y\inv$, we obtain
\eqnref{almost}.
\end{proof}

\begin{lemma}
\lemlabel{dagger}
For all $x,y\in Q$,
\begin{equation}
\eqnlabel{dagger}
L_x\inv D_y L_x = L_{xy}\inv S_{(xy)\ld x} L_{xy}\,.
\end{equation}
\end{lemma}

\begin{proof}
We compute
\begin{alignat*}{2}
L_x\inv D_y L_x
&\ = L_x\inv L_y\inv S_y L_y L_x
&&\ = L_x\inv L_y\inv S_y L_{y,x} L_{xy} \\
&\by{AL} L_x\inv L_y\inv L_{y,x} S_{y L_{y,x}} L_{xy}
&&\by{fix-y} L_{xy}\inv S_{(xy)\ld x} L_{xy}\,,
\end{alignat*}
where we have also used the assumption that $Q$ has exponent $2$
in the last step.
\end{proof}

Finally, we have enough for the main result of this section.

\begin{proof}[Proof of Theorem \thmref{abelian-group}]
By commutativity of $\nlp$
(Proposition \prpref{newop-comm-loop}) and $x\nlp x = x^2 = 1$
for all $x\in Q$, all that is needed is to show that $\nlp$ is
associative. First, apply \eqnref{dagger} to both sides of \eqnref{almost}
to obtain
$L_{xy}\inv S_{(xy)\ld x} L_{xy} = L_{yx}\inv S_{(yx)\ld y} L_{yx} D_{xy}$,
or $S_{(xy)\ld x} = S_{(yx)\ld y} L_{yx} D_{xy} L_{xy}\inv = S_{(yx)\ld y} S_{xy}$.
Replace $x$ with $y\ld x$ to get
$S_{x\ld (y\ld x)} = S_{x\ld y} S_x$. Replace $y$ with $xy$ to obtain
$S_{x\ld (xy\ld x)} = S_y S_x$, or $S_{x\nlp y} = S_y S_x$. This is
precisely associativity of $\nlp$: applying both sides to $z$, we have
$(x\nlp y)\nlp z = x\nlp (y\nlp z)$ for all $x,y,z\in Q$. This completes
the proof.
\end{proof}

\section{$p$-loops}
\seclabel{p-loops}

For a finite, power-associative loop $Q$, there are at least two reasonable ways to
define what it means for $Q$ to be a $p$-\emph{loop}: either every element of $Q$
has order a power of $p$, or $|Q|$ is a power of $p$. Fortunately,
these two notions are equivalent for groups, Moufang loops, and, as we are
about to show, for commutative A-loops.

\begin{theorem}
\thmlabel{p-loop}
Let $Q$ be a finite commutative A-loop and let $p$ be a prime.
Then $|Q|$ is a power of $p$ if and only if every element of $Q$
has order a power of $p$.
\end{theorem}

\begin{proof}
Assume first that $p$ is odd. If $|Q|$ is a power of $p$, then by Proposition
\prpref{lagrange}, every element of $Q$ has order a power of $p$. Conversely,
if $|Q|$ is divisible by an odd prime $q$, then by Proposition
\prpref{HSC}(iii), $Q$ contains an element of order $q$. Thus if every element
of $Q$ has order a power of $p$, $|Q|$ must be a power of $p$.

Now assume that $p = 2$ and that $|Q|$ is a power of $2$.
Since $Q = K(Q)\times H(Q)$ (Theorem \thmref{Decomposition})
and $|K(Q)|$ is odd (Theorem \thmref{IteratedSquares}),
we must have  $K(Q) = \sbl{1}$, and so $Q = K(Q)$, that is,
every element of $Q$ has order a power of $2$.

For the converse, assume that $Q$ is a smallest commutative A-loop of
exponent a power of $2$ such that $|Q|$ is not a power of $2$. Consider
the normal subloop $1 < H_1 = \setof{x\in Q}{x^2 = 1}$, \emph{cf.}
Theorem \thmref{Power2Exp}. Then $|H_1|$ is a power of $2$ by Corollary
\corref{exp2-order}. If $H_1 = Q$, we have reached a contradiction.
If $H_1 < Q$ then $|Q/H_1|$ is a power of $2$ by minimality, and so
$|Q| = |H_1|\cdot |Q/H_1|$ is a power of $2$, a contradiction.
\end{proof}

Unlike in the case of abelian groups, for a finite commutative
A-loop $Q$, the normal subloop $K(Q)$ does not necessarily decompose
as a direct product of $p$-loops. For example, Dr\'apal \cite{Drapal}
constructed a commutative A-loop of order $15$ that is not a direct
product of a $3$-loop and a $5$-loop.

\begin{theorem}[Lagrange and Cauchy Theorems]
\thmlabel{LHSC} Let $Q$ be a finite commutative A-loop. Then:
\begin{enumerate}[label=\roman*\emph{)}]
\item\quad If $x\in A\le Q$ then both $|x|$ and $|A|$ divide $|Q|$.
\item\quad If a prime $p$ divides $|Q|$ then $Q$ has an element of order
    $p$.
\end{enumerate}
\end{theorem}
\begin{proof}
Combine Theorems \thmref{Decomposition}, \thmref{p-loop} and Propositions
\prpref{lagrange}, \prpref{HSC}.
\end{proof}

\section{Open Problems}
\seclabel{open}

We conclude this paper with some open problems.

\begin{problem}
\prblabel{simple}
Does there exist a nonassociative, finite simple commutative
A-loop?
\end{problem}

By Proposition \prpref{simple} and Corollary \corref{exp2-order},
such a loop would have exponent $2$ and order a power of $2$.
To get some insight into the problem, more constructions of
commutative A-loops which are $2$-loops are needed; see
\cite{JKV2}.

Recall that the \emph{center} of a loop $Q$ is the set of all elements $a$
satisfying $a\cdot xy = x\cdot ay = xa\cdot y$ for all $x,y$. In groups and
Moufang loops, the center of a $p$-loop is always nontrivial, and thus such
loops are centrally nilpotent.

\begin{problem}
\prblabel{center}
Let $p$ be an odd prime. Does there exist a finite commutative
A-loop of order a power of $p$ with trivial center?
\end{problem}

By a classic result of Albert \cite{Albert}, it would be
sufficient to show that $\Mlt(Q)$ is a $p$-group.

The restriction to odd $p$ is necessary.
There exist commutative A-loops of exponent $2$ of all
orders $2^n$, $n \geq 4$ with trivial center \cite{JKV2}.

For a set $\pi$ of primes, a positive integer $n$ is a $\pi$-\emph{number}
if $n = 1$ or if $n$ is a product of primes in $\pi$. For each positive
integer $n$, let $n_{\pi}$ denote the largest $\pi$-number dividing $n$.
A subloop $K$ of a finite, power-associative
loop $Q$ is a \emph{Hall} $\pi$-\emph{subloop} if $|K| = |Q|_{\pi}$. In case
$\pi = \{p\}$, we say that $K$ is a \emph{Sylow} $p$-\emph{subloop} of $Q$.

\begin{problem}
\prblabel{hall}
Let $Q$ be a commutative A-loop.
\begin{enumerate}[label=\roman*\emph{)}]
\item\quad For each set $\pi$ of primes, does $Q$ have a Hall $\pi$-subloop?
\item\quad For each prime $p$, does $Q$ have a Sylow $p$-subloop?
\end{enumerate}
\end{problem}

Sylow $2$-subloops certainly exist by Theorems \thmref{Decomposition} and
\thmref{p-loop}. In both the Hall and Sylow cases, the problem reduces to
considering commutative A-loops of odd order. Hall and Sylow subloops of the
associated Bruck loop $(Q,\circ)$ exist \cite{Glauberman1}, so the question is whether
or not these are also subloops of $Q$ itself.

\end{document}